\documentclass[12pt,leqno]{article}
\usepackage[official]{eurosym}
\usepackage{fullpage,amsmath,theorem}
\usepackage{epic,eepic}
\usepackage{amssymb,amsmath}
\usepackage{cases}
\def\llvdash{{\|\hskip-2pt \raise 3pt\hbox{\vrule
height 0.25pt width 0.4cm}}}

\def\l{{\langle}}
\def\r{{\rangle}}

\def\calP{\mathcal P}

\def\calK{\mathcal K}

\def\upr{\upharpoonright}

\def\oa{{\overline A^{\,\lower 7pt_{\hbox{$\scriptstyle\bet}}
\hbox{$\scriptstyle 0\tau$}}}}

\def\bet{\beta}

\def\llvdash{{\|\hskip-2pt \raise 3pt\hbox{\vrule height 0.25pt
width 0.4cm}}}

\newtheorem{theorem}{Theorem}[section]
\newtheorem{lemma}[theorem]{Lemma}
\newtheorem{corollary}[theorem]{Corollary}
\newtheorem{proposition}[theorem]{Proposition}
{\theorembodyfont{\rmfamily}
}
{\theorembodyfont{\rmfamily}
\newtheorem{remark}[theorem]{Remark}}
{\theorembodyfont{\rmfamily}
\newtheorem{claim}{Claim}}
{\theorembodyfont{\rmfamily}
}
{\theorembodyfont{\rmfamily}
}
\DeclareMathOperator{\pcf}{pcf} 
\DeclareMathOperator{\rng}{rng} \DeclareMathOperator{\cof}{cof}

\newcommand{\pr}{\medskip\noindent\textit{Proof}. }

\newcommand{\lusim}[1]{\smash{\underset{\raisebox{1.2pt}[0cm][0cm]{$\sim$}}
{{#1}}}}
\def\tcf{{\rm tcf}}
\def\bd{{\rm bd}}
\def\pp{{\rm pp}}

\def\cf{{\rm cf}}
\def\rng{{\rm rng}}

\def\complete{{\rm complete}}
\def\Reg{{\rm Reg}}

\def\llvdash{{\|\hskip-2pt \raise 3pt\hbox{\vrule
height 0.25pt width 0.15cm}}}

\def\Vdashbks{\hbox{$\Vdash\!\!\!\!{\raise2pt\hbox
{$\scriptscriptstyle\backslash$}}$}}

\begin{document}

\title{Applications of pcf for mild large cardinals to elementary embeddings. }
\baselineskip=18pt
\author{ Moti Gitik\footnote{ We are grateful to Menachem Magidor
for his comments. Gitik was partially supported by ISF grant 234/08} and Saharon Shelah\footnote{ Shelah was partially supported by ISF grant 1053/11. This is paper 1013 on Shelah's   publication list.}
 }

\date{}
\maketitle

\begin{abstract}
The following
pcf results are proved:

\emph{1. Assume that $\kappa>\aleph_0$ is a weakly compact cardinal.
Let $\mu>2^\kappa$ be a singular cardinal of cofinality $\kappa$.
Then for every regular $\lambda<\pp^+_{\Gamma(\kappa)}(\mu)$
there is an increasing sequence $\l \lambda_i \mid i<\kappa \r$ of regular cardinals converging to $\mu$
such that
$\lambda= \tcf(\prod_{i<\kappa} \lambda_i, <_{J^{bd}_\kappa})$.}

\emph{2.
Let $\mu$ be a strong limit cardinal and $\theta$  a cardinal above $\mu$.
Suppose that at least one of them has an uncountable cofinality.
Then   there is $\sigma_*<\mu$
such that for every $\chi<\theta$ the following holds:
$$ \theta> \sup\{ \sup \pcf_{\sigma_*-\complete}(\frak{a}) \mid
\frak{a}\subseteq \Reg \cap (\mu^+,\chi) \text{ and } |\frak{a}|<\mu \}.$$}

As an application we show
that: \\\emph{if $\kappa$ is a measurable cardinal and $j:V \to M$ is the
elementary embedding by a $\kappa$--complete ultrafilter over a
measurable cardinal $\kappa$, then for every $\tau$ the following
holds:
\begin{enumerate}
  \item   if $j(\tau)$ is a cardinal then
$j(\tau)=\tau$;
  \item $|j(\tau)|=|j(j(\tau))|$;
  \item for any $\kappa$--complete ultrafilter $W$ on $\kappa$,\quad
$|j(\tau)|=|j_W(\tau)|$.
\end{enumerate}}
The first two items provide affirmative answers to   questions from
\cite{G-Sh} and the third to
  a question of D. Fremlin.
\end{abstract}

\section{Introduction}

We address here the following question:

Suppose $\kappa$ is a measurable cardinal, $U$ a $\kappa$--complete non-trivial ultrafilter
over $\kappa$ and $j:V \to M$ the corresponding elementary embedding.
Can one characterize cardinals moved by $j$?

There are trivial answers. For example:
\\$\tau$ is moved by $j$ iff $\cof(\tau)=\kappa$ or there is some $\delta<\tau$ with $j(\delta)\geq \tau$.
\\Also, assuming GCH, it is not hard to find a characterization in terms not mentioning
$j$.

However, it turns out that an answer is possible in terms not
mentioning $j$ already in ZFC (Theorem \ref{thm0-10}):

\emph{Let $\tau$ be a cardinal.
 Then either
\begin{enumerate}
\item $\tau<\kappa$ and then $j(\tau)=\tau$,\\ or
\item $\kappa\leq \tau \leq 2^\kappa$ and then $j(\tau)>\tau$,
 $2^\kappa<j(\tau)<(2^\kappa)^+$,\\ or
 \item $\tau\geq (2^\kappa)^+$ and then $j(\tau)>\tau$ iff
there is a singular cardinal ${\mu}\leq \tau$  of cofinality $\kappa$ above $2^\kappa$ such that
$\pp_{\Gamma(\kappa)}({\mu})\geq\tau$, and if $\tau^*$ denotes the least such $\mu$, then\\
$\tau\leq\pp_{\Gamma(\kappa)}(\tau^*)<j(\tau)<\pp_{\Gamma(\kappa)}(\tau^*)^+$.
\end{enumerate}
}

Straightforward  conclusions of this result provide
affirmative answers to questions mentioned in the abstract.

A crucial tool here is PCF--theory and specially Revisited GCH Theorem \cite{Sh:460} Sh460.
\\A new result involving weakly compact cardinal is obtained (Theorem \ref{pcf-main}):

\emph{Assume that $\kappa>\aleph_0$ is a weakly compact cardinal.
Let $\mu>2^\kappa$ be a singular cardinal of cofinality $\kappa$.
Then for every regular $\lambda<\pp^+_{\Gamma(\kappa)}(\mu)$
there is an increasing sequence $\l \lambda_i \mid i<\kappa \r$ of regular cardinals converging to $\mu$
such that
$\lambda= \tcf(\prod_{i<\kappa} \lambda_i, <_{J^{bd}_\kappa})$.
}

Also a bit sharper version of \cite{Sh:460} Sh460, 2.1 for uncountable cofinality is proved (Theorem \ref{thm2}):

\emph{
Let $\mu$ be a strong limit cardinal and $\theta$  a cardinal above $\mu$.
Suppose that at least one of them has an uncountable cofinality.
Then   there is $\sigma_*<\mu$
such that for every $\chi<\theta$ the following holds:
$$ \theta> \sup\{ \sup \pcf_{\sigma_*-\complete}(\frak{a}) \mid
\frak{a}\subseteq \Reg \cap (\mu^+,\chi) \text{ and } |\frak{a}|<\mu \}.$$
}

The first author proved a version of \ref{thm0-10} assuming certain weak form of the Shelah Weak Hypothesis (SWH)\footnote{
Consistency of negations of SWH is widely open except very few instances.}
and using \cite{Sh:g} Sh371.
Then the second author was able to show that the actual assumption used holds in ZFC.
All PCF results of the paper are due solely to him.

Let us recall the definitions of few basic notions of PCF theory that will be used here.

Let $\frak{a}$ be a set of regular cardinals above $|\frak{a}|$.

$$\pcf(\frak{a})=\{\tcf((\prod \frak{a}, <_J))\mid J \text{ is an  ideal on }
\frak{a}$$$$ \text{ and }  (\prod \frak{a}, <_J) \text{ has true cofinality }\}.$$

Let $\rho$ a cardinal.

$$\pcf_{\rho-\complete}(\frak{a})=\{\tcf((\prod \frak{a}, <_J))\mid J \text{ is a } \rho-\text{complete ideal on }
\frak{a}$$$$ \text{ and }  (\prod \frak{a}, <_J) \text{ has true cofinality }\}.$$

Let $\eta$ be a cardinal.

$$J_{<\eta}[\frak{a}]=\{\frak{b}\subseteq \frak{a} \mid \text{ for every ultrafilter } D \text{ on }
\frak{b}, \cf(\prod \frak{b}, <_D)<\lambda\}.$$

Let $\lambda$ be a singular cardinal.
$$\pp_{\Gamma(\kappa)}(\lambda)=\pp_{\Gamma(\kappa^+,\kappa)}(\lambda)=\sup\{tcf((\prod \frak{a}, <_J))\mid
\frak{a} \text{ is a set of } \kappa \text{ regular cardinals unbounded in } \lambda,$$
$$J \text{ is a } \kappa-\text{complete ideal on } \frak{a} \text{ which includes } J_{\frak{a}}^{\bd}
\text{ and } (\prod \frak{a}, <_J) \text{ has true cofinality }\}.$$

$$\pp_{\Gamma(\kappa)}^+(\lambda) \text{ denotes the first regular without such representation.
}\footnote{ Note that $\pp_{\Gamma(\kappa)}^+(\lambda)\leq
(\pp_{\Gamma(\kappa)}(\lambda))^+$ and it is open if
$\pp_{\Gamma(\kappa)}^+(\lambda)< (\pp_{\Gamma(\kappa)}(\lambda))^+$
can ever occur (see \cite{Sh:g},Sh355, p.41.)}$$

\section{PCF results.}

\begin{theorem}\label{pcf-main}
Assume that $\kappa>\aleph_0$ is a weakly compact cardinal.
Let $\mu>2^\kappa$ be a singular cardinal of cofinality $\kappa$.
Then for every regular $\lambda<\pp^+_{\Gamma(\kappa)}(\mu)$
there is an increasing sequence $\l \lambda_i \mid i<\kappa \r$ of regular cardinals converging to $\mu$
such that
$\lambda= \tcf(\prod_{i<\kappa} \lambda_i, <_{J^{bd}_\kappa})$.

\end{theorem}

\begin{remark}\label{rem-pcf-main}
It is possible to remove the assumption $\mu>2^\kappa$.
Just \cite{Sh:430}(Sh430) $\S$ 6, 6.7A should be used to find the pcf-generators in the proof below.
See also 6.3 of Abraham -Magidor handbook article \cite{A-M}.
\end{remark}

\pr
By No Hole Theorem (2.3, p.57 \cite{Sh:g}),  there are
a $\kappa$--complete ideal $I_1$ on $\kappa$ and a sequence of regular cardinals $\vec{\lambda}^1=\l \lambda^1_i\mid i<\kappa \r$
with $\mu=\lim_{I_1}\vec{\lambda}^1$ such that \\
$\lambda= \tcf(\prod_{i<\kappa} \lambda_i^1, <_{I_1})$.

Denote the set $\{ \lambda^1_i\mid i<\kappa \}$ by $\frak{a}^1$. Let
$\frak{a}^2=\pcf(\frak{a}^1)$.
\\Without loss of generality assume that $\lambda=\max \pcf(\frak{a}^1)$.
Note that  by \cite{Sh:g} the following holds:
\begin{enumerate}
\item
$\frak{a}^1\subseteq \frak{a}^2\subseteq Reg \setminus \kappa^{++}$,
\item $\pcf(\frak{a}^2)= \frak{a}^2$,
\item $|\pcf(\frak{a}^2)|\leq 2^\kappa$.
\end{enumerate}

By \cite{Sh:g}([Sh345a, 3.6, 3.8(3)) there is a smooth and closed generating sequence for $\frak{a}^1$
(here we use $2^\kappa<\mu$)
which means a sequence $\l \frak{b}_\theta \mid \theta\in \frak{a}^2\r$ such that
\begin{enumerate}
  \item $\theta\in \frak{b}_\theta \subseteq \frak{a}^2$,
  \item $\theta \not \in \pcf(\frak{a}^2 \setminus \frak{b}_\theta)$,
  \item $\frak{b}_\theta=\pcf(\frak{b}_\theta)$,
  \item $\theta_1 \in \frak{b}_{\theta_2}$ implies $\frak{b}_{\theta_1} \subseteq \frak{b}_{\theta_2}$,
  \item $\theta=\max \pcf(\frak{b}_\theta)$.
\end{enumerate}

Then by \cite{Sh:g}[Sh345a,3.2(5)]:

(*)$_1$: \quad if $\frak{c}\subseteq \frak{a}^2$, then for some finite $\frak{d}\subseteq \pcf(\frak{c})$ we have
$\frak{c}\subseteq \pcf(\frak{c})\subseteq \bigcup\{\frak{b}_\theta \mid \theta \in \frak{d}\}$.

The next claim is a consequence of \cite{Sh:460}(Sh460, 2.1):

\begin{claim}\label{clm1}
There is $\sigma_*<\kappa$ such that for every $\frak{a} \subset Reg \cap (\kappa^+,\mu)$ of cardinality less than $\kappa$
there is a sequence $\l \frak{a}_\alpha \mid \alpha<\sigma_* \r$ such that
\begin{enumerate}
  \item $\frak{a} = \bigcup_{\alpha<\sigma_*} \frak{a}_\alpha$,
  \item $\max \pcf(\frak{a}_\alpha)< \mu$, for every $\alpha<\sigma_*$.
\end{enumerate}
\end{claim}
\pr
The cardinal $\kappa$ is a strong limit, so we can apply \cite{Sh:460}(Sh460, 2.1)
to $\kappa$ and $\mu$. Hence there is $\sigma_*<\kappa$ such that for every $\frak{a} \subset Reg \cap (\kappa^+,\mu)$ of cardinality less than $\kappa$ we have $\pcf_{\sigma_*^+-\complete}(a)\subseteq \mu$.
This means that the $\sigma_*^+$--complete ideal generated by $J_{<\mu}(\frak{a})$ is everything, i.e. $\calP(\frak{a})$.
See 8.5 of \cite{A-M} for the detailed argument.
So there are
$ \frak{a}_\alpha$'s  in $J_{<\mu}(\frak{a})$, for $\alpha<\sigma_* $ such that
 $\frak{a} = \bigcup_{\alpha<\sigma_*} \frak{a}_\alpha$. But then also
   $\max \pcf(\frak{a}_\alpha)< \mu$, for every $\alpha<\sigma_*$.
\\
$\square$ of the claim.

Let $\sigma_*<\kappa$ be given by the claim.
Let $i<\kappa$.
Apply the  claim to the set $\frak{a}^1_i:=\{\lambda^1_j \mid j<i \}$.
So there is a sequence $\l \frak{a}_{i\alpha} \mid \alpha<\sigma_* \r$ such that
\begin{enumerate}
  \item $\frak{a}^1_i = \bigcup_{\alpha<\sigma_*} \frak{a}_{i\alpha}$,
  \item $\max \pcf(\frak{a}_{i\alpha})< \mu$, for every $\alpha<\sigma_*$.
\end{enumerate}
Now, by (*)$_1$, for every $\alpha<\sigma_*$,
$$\pcf(\frak{a}_{i\alpha})\subseteq \bigcup\{\frak{b}_\theta \mid \theta \in \frak{d}_{i\alpha}\},$$
for some finite $\frak{d}_{i\alpha}\subseteq \pcf(\frak{a}_{i\alpha})$.
\\Set $\frak{d}_i=\bigcup_{\alpha<\sigma_*} \frak{d}_{i\alpha}$.
Then $\frak{d}_i$ is a subset of $\mu$ of cardinality $\leq \sigma_*$.
In addition we have
$\frak{d}_i\subseteq \pcf(\frak{a}_{i}^1)$ and $\frak{a}_{i}^1 \subseteq \bigcup\{\frak{b}_\theta \mid \theta \in \frak{d}_{i}\}$.
\\Let $\l \theta_{i,\epsilon} \mid \epsilon<\sigma_* \r$ be a listing of $\frak{d}_{i}$.

\begin{claim}\label{clm2}
There are a function ${g}$ and $\vec{u}=\l u_\epsilon \mid \epsilon<\sigma_* \r$
such that
\begin{enumerate}
  \item $g:\kappa \to \kappa$ is increasing,
  \item $\xi\leq g(\xi)$, for every $\xi<\kappa$,
  \item $\kappa= \bigcup_{\epsilon<\sigma_*}u_\epsilon$,
  \item for any $\epsilon<\sigma_*$ and $\xi<\eta<\kappa$ the following holds:
  $$\lambda^1_{\xi}\in \frak{b}_{\theta_{g(\eta),\epsilon}} \text{ iff } \xi \in u_\epsilon.$$

\end{enumerate}

\end{claim}
\pr
Here is the place to use the weak compactness of $\kappa$.

We will define  a $\kappa$--tree $T$ and then will use its $\kappa$--branch.

Fix $\eta<\kappa$.
Let $P \subseteq \sigma_*\times\eta $.
Define a set
$$A_P:=\{\alpha\in (\eta,\kappa) \mid \forall \xi<\eta \forall \epsilon<\sigma_*
(\l\epsilon,\xi\r \in P \Leftrightarrow\lambda^1_\xi \in \frak{b}_{\theta_{\alpha,\epsilon}})\}.$$
Note that always there is $P \subseteq \sigma_*\times\eta $ with $|A_P|=\kappa$.
Just $|\calP(\sigma_*\times\eta)|<\kappa$, so the function
$$\alpha\longmapsto \l \l\epsilon,\xi\r \mid \epsilon<\sigma_*, \xi<\eta \text{ and }\lambda^1_\xi \in \frak{b}_{\theta_{\alpha,\epsilon}}\r$$
is constant on a set of cardinality $\kappa$.
\\Also for such $P$ we will have $\rng(P)=\eta$, i.e. for every $\xi<\eta$ there is $\epsilon<\sigma_*$
(which may be not unique)
such that $(\epsilon,\xi) \in P$.
Thus pick  $\alpha \in A_P$. Then $\alpha>\eta>\xi$ and
$\frak{a}_{\alpha}^1 \subseteq \bigcup\{\frak{b}_\theta \mid \theta \in \frak{d}_{\alpha}\}$.
Clearly $\lambda^1_\xi$ appears in $\frak{a}_{\alpha}^1=\{\lambda^1_\nu \mid \nu<\alpha\}$.
Hence there is $\epsilon<\sigma_*$ such that $\lambda^1_\xi \in b_{\theta_{\alpha,\epsilon}}$, and so
$(\epsilon,\xi) \in P$.

Let $$T:=\{P \mid \exists \eta<\kappa (P \subseteq \sigma_*\times\eta  \text{ and } |A_P|=\kappa) \}.$$
If $P \subseteq \sigma_*\times\eta, P' \subseteq \sigma_*\times\eta'$ are both in  $ T $  then set
$P <_T P' $ iff
\begin{itemize}
  \item $\eta<\eta'$,
  \item $P'\cap (\sigma_*\times\eta)=P$.
\end{itemize}

 Then $\l T, <_T \r$ is a $\kappa$--tree.
Let $X \subseteq \sigma_* \times\kappa   $ be a $\kappa$--branch.
Define now an increasing function $g:\kappa \to \kappa$.
Set $g(\eta)=\min(A_{X \cap (\sigma_*\times \eta)}\setminus \sup\{g(\eta')\mid \eta'<\eta\})$.
\\Let now $\epsilon<\sigma_*$. Define $u_\epsilon$ as follows:
$$\xi\in u_\epsilon \text{ iff for some } \eta>\xi \text{ and some (every)} \alpha \in A_{X \cap (\sigma_*\times \eta)},
\lambda^1_\xi \in \frak{b}_{\theta_{\alpha,\epsilon}}.$$

Then
for any $\epsilon<\sigma_*$ and $\xi<\eta<\kappa$ the following holds:
  $$\lambda^1_{\xi}\in \frak{b}_{\theta_{g(\eta),\epsilon}} \text{ iff } \xi \in u_\epsilon.$$
Finally $|X|=\kappa$ implies that for every $\xi<\kappa$ there is $\epsilon<\sigma_*$ with $\xi\in u_\epsilon$.
Thus let $\xi<\kappa$. Pick some $\eta, \xi<\eta<\kappa$. Consider $X \cap (\sigma_*\times \eta)$.
Then, as was observed above, there are $\alpha \in A_{ X \cap (\sigma_*\times \eta)}$ and $\epsilon<\sigma_*$
such that $\lambda^1_\xi \in b_{\theta_{\alpha,\epsilon}}$. Hence $\xi\in u_\epsilon$.
\\
$\square$ of the claim.

\begin{claim}\label{clm3}
Suppose that $u_\epsilon \in I_1^+$, for some $\epsilon<\sigma_*$.
Then $|u_\epsilon|=\kappa$ and the quasi order $\prod_{i \in u_\epsilon} (\theta_{g(i),\epsilon},<_{J^{\bd}_{u_\epsilon}})$
has true cofinality $\lambda$.

\end{claim}
\pr
 $\kappa$--completeness of $I_1$ implies that $|u_\epsilon|=\kappa$, since clearly $\{\xi\} \in I_1$, for every $\xi<\kappa$.
\\Suppose now that
the quasi order $\prod_{i \in u_\epsilon} (\theta_{g(i),\epsilon},<_{J^{\bd}_{u_\epsilon}})$ does not have a true cofinality or it
has true cofinality $\not =\lambda$.
Recall that $\lambda=\max \pcf(\frak{a}_1)$.
So by \cite{Sh:g}(Sh345a) there is an unbounded subset $v $ of $u$ such that
$\prod_{i \in v} (\theta_{g(i),\epsilon},<_{J^{\bd}_{v}})$ has a true cofinality
$\lambda_*<\lambda$.
We can take $\lambda_*$ to be just the least $\delta$ such that an unbounded subset of $u_\epsilon$ appears in $J_{\leq \delta}[u_\epsilon]$.
Without loss of generality we can assume that $\lambda_*=\max \pcf(\{\theta_{g(i),\epsilon}\mid i \in v\})$.
We have
$\lambda_*\in \pcf(\{\theta_{g(i),\epsilon}\mid i \in v\})\subseteq \pcf(\frak{a}_1)=\frak{a}_2$.
Set $v_1:=\{i \in v \mid \theta_{g(i),\epsilon}\in \frak{b}_{\lambda_*}\}$.
Then $v_1$ is unbounded in $v$. By smoothness of the generators, \quad $i \in v_1$ implies $\frak{b}_{\theta_{g(i),\epsilon}}\subseteq \frak{b}_{\lambda_*}$.
Then
$$i \in v_1 \text{ and } \xi\in u_\epsilon\cap i \text{ imply } \lambda^1_\xi \in \frak{b}_{\lambda_*}.$$
But $v_1$ is unbounded in $\kappa$, hence for every $\xi \in u_\epsilon$ there is $i \in v_1, i>\xi$.
So, $\{\lambda^1_\xi\mid \xi \in u_\epsilon\} \subseteq \frak{b}_{\lambda_*}$.
By the closure of the generators,  $\pcf(\frak{b}_{\lambda_*})=\frak{b}_{\lambda_*}$.
Hence $\pcf(\{\lambda^1_\xi\mid \xi \in u_\epsilon\} )\subseteq \frak{b}_{\lambda_*}$.
This impossible since  $u_\epsilon \in I_1^+$ and so $\lambda \in \pcf(\{\lambda^1_\xi\mid \xi \in u_\epsilon\} )$,
but $\lambda_*<\lambda$. Contradiction.
\\
$\square$ of the claim.

\begin{claim}\label{clm4}
There is $\epsilon<\sigma_*$ such that $u_\epsilon \in I_1^+$ and $\mu=\lim_{J^{\bd}_\kappa+(\kappa\setminus u_\epsilon)}\l
\theta_{g(i),\epsilon} \mid i<\kappa\r$.

\end{claim}
\pr
Suppose otherwise.
Set $s:=\{\epsilon<\sigma_* \mid u_\epsilon \in I_1^+\}$.
Then for every $\epsilon\in s $ there is $v_\epsilon $ an unbounded subset of $\kappa$ such that
$\theta_\epsilon^*:=\sup\{\theta_{g(i),\epsilon}\mid i \in v_\epsilon\}$ is below $\mu$.
Set \\$\theta_*:=\sup\{\theta^*_\epsilon \mid \epsilon\in s\}$.
Then $\theta_*<\mu$, since $\cof(\mu)=\kappa>\sigma_*$.
\\
Set $w_1:=\bigcup\{u_\epsilon \mid \epsilon \in \sigma_*\setminus s\}$.
Then $w_1 \in I_1$ as a union of less than $\kappa$ of its members.
Also the set $w_2:=\{i<\kappa \mid \lambda^1_i \leq\theta_* \} $ belongs to $I_1$ because
$\mu=\lim_{I_1}\{\lambda^1_i \mid i<\kappa\}$.
Hence $w:=w_1\cup w_2 \in I_1$.
\\Let $\xi \in \kappa \setminus w$.
Then
$$\lambda_\xi^1 \in \{\lambda_\rho^1 \mid \rho<\xi+1\} \subseteq \bigcup\{\frak{b}_{\theta_{g(\xi+1),\epsilon}} \mid \epsilon<\sigma_*\}.$$
Hence for some $\epsilon<\sigma_*$, \quad $\lambda_\xi^1\in \frak{b}_{\theta_{g(\xi+1),\epsilon}}$.
Then $\xi \in u_\epsilon$.
Now, $\xi \not \in w$ and so $\xi \not \in w_1$. Hence $\epsilon \in s$.
Pick some $\tau \in v_\epsilon, \tau>\xi$. Then $\lambda^1_\xi \in \frak{b}_{\theta_{{g(\tau),\epsilon}}}$, since $\xi \in u_\epsilon$.
Then $$\lambda^1_\xi \leq \max( \frak{b}_{\theta_{{g(\tau),\epsilon}}})=\theta_{{g(\tau),\epsilon}}
\leq \theta^*_\epsilon\leq \theta_*.$$
But then $\xi \in w_2$. Contradiction.
\\
$\square$ of the claim.
\\
$\square$

\begin{proposition}\label{lem-x}
Let $\frak{a}$ be a set of regular cardinals with $\min(\frak{a})>2^{|\frak{a}|}$.
Let $\sigma<\theta\leq |\frak{a}|$.
Suppose that $\lambda\in \pcf_{\sigma-\complete}(\frak{a})$, $\mu<\lambda$ and
$\pcf_{\theta-\complete}(\frak{a})\subseteq \mu$.
Then there is $\frak{c} \subseteq \pcf_{\theta-\complete}(\frak{a})$ such that
$|\frak{c}|<\theta$,
   $\frak{c}\subseteq \mu$ and
   $\lambda \in \pcf_{\sigma-\complete}(\frak{c})$.

\end{proposition}

\begin{remark}
It is possible to replace the assumption $\min(\frak{a})>2^{|\frak{a}|}$ by $\min(\frak{a})>|\frak{a}|$
using \cite{Sh:430}(Sh430) $\S$ 6, 6.7A  in order to find the pcf-generators used in the proof.

\end{remark}
\pr
Let
$\l \frak{b}_\xi \mid \xi \in \pcf(\frak{a})\r$ be a set of generators as in Theorem \ref{pcf-main}.
We have $\lambda\in \pcf_{\sigma-\complete}(\frak{a})\subseteq \pcf(\frak{a})$, hence
$\frak{b}_\lambda$ is defined and $\max \pcf(\frak{b}_\lambda)=\lambda \in \pcf_{\sigma-\complete}(\frak{a})\subseteq \pcf(\frak{a})$.
\\By \cite{Sh:430}, 6.7F(1), there is $\frak{c}\subseteq \pcf_{\theta-\complete}(\frak{a}\cap \frak{b}_\lambda)\subseteq \mu$ of cardinality $<\theta$ such that \\
$ \frak{b}_\lambda\cap \frak{a}\subseteq \bigcup\{\frak{b}_\xi \mid \xi \in \frak{c}\}.$
Then, by smoothness, $\xi \in \frak{c}\Rightarrow \frak{b}_\xi \subseteq \frak{b}_\lambda$.
Also $\pcf(\frak{c})\subseteq \pcf(\frak{b}_\lambda)=\frak{b}_\lambda$. Hence $\max \pcf(\frak{c})\leq \lambda$.

Now, if $\lambda\in \pcf_{\sigma-\complete}(\frak{c})$, then we are done.
Suppose otherwise.
Then there are $j(*)<\sigma$ and $\theta_j\in \lambda\cap \pcf_{\sigma-\complete}(\frak{c})$, for every $j<j(*)$,
such that
$\frak{c}\subseteq \bigcup\{\frak{b}_{\theta_j} \mid j<j(*)\}$.
So if $\eta \in \frak{b}_\lambda \cap \frak{a}$, then for some $\chi \in \frak{c}$
we have $\eta \in \frak{b}_\chi$, as $ \frak{b}_\lambda \cap \frak{a}\subseteq \bigcup\{\frak{b}_\xi \mid \xi \in \frak{c}\}.$
Hence for some $j<j(*)$, \quad $\chi \in \frak{b}_{\theta_j}$, and so
 $\frak{b}_\chi \subseteq \frak{b}_{\theta_j}$ and  $\eta \in \frak{b}_{\theta_j}$.
\\Then $\frak{b}_\lambda \cap \frak{a}\subseteq \bigcup_{j<j(*)} \frak{b}_{\theta_j}$.
Recall that $j(*)<\sigma$ and $\theta_j<\lambda$, for every $j<j(*)$.
\\Note that $\lambda\in \pcf_{\sigma-\complete}(\frak{a})$ implies that
$\lambda\in \pcf_{\sigma-\complete}(\frak{b}_\lambda\cap \frak{a})$, see for example 4.14 of \cite{A-M}.
So there is a $\sigma$--complete ideal $J$
on $\frak{b}_\lambda\cap \frak{a}$ such that\\ $\lambda=\tcf(\prod(\frak{b}_\lambda\cap \frak{a}),<_J)$.
Then for some $j<j(*)$, \quad $\frak{b}_{\theta_j} \in J^+$ which is impossible since
$\max \pcf(\frak{b}_{\theta_j})=\theta_j<\lambda$. Contradiction.
\\
$\square$

The next result follows from 2.1 of \cite{Sh:460} Sh460.

\begin{theorem}\label{thm2}
Let $\mu$ be a strong limit cardinal and $\theta$  a cardinal above $\mu$.
Suppose that at least one of them has an uncountable cofinality.
Then   there is $\sigma_*<\mu$
such that for every $\chi<\theta$ the following holds:
$$ \theta> \sup\{ \sup \pcf_{\sigma_*-\complete}(\frak{a}) \mid
\frak{a}\subseteq \Reg \cap (\mu^+,\chi) \text{ and } |\frak{a}|<\mu \}.$$

\end{theorem}

\pr
Assume first that $\cof(\mu)\not = \cof(\theta)$.
Suppose on contrary that
$$\forall \mu^*<\mu \exists \chi<\theta
(\theta\leq \sup\{ \sup \pcf_{\mu^*-\complete}(\frak{a}) \mid
\frak{a}\subseteq \Reg \cap (\mu^+,\chi) \text{ and } |\frak{a}|<\mu \}).$$
If $\cof(\theta)<\cof(\mu)$, then there will be $\chi<\theta$ such that for every $\mu^*<\mu$
$$\theta\leq \sup\{ \sup \pcf_{\mu^*-\complete}(\frak{a}) \mid
\frak{a}\subseteq \Reg \cap (\mu^+,\chi) \text{ and } |\frak{a}|<\mu \}.$$
But this is impossible by 2.1 of \cite{Sh:460} applied to $\mu$ and $\chi$.
\\If $\cof(\theta)>\cof(\mu)$, then still
there will be $\chi<\theta$ such that for every $\mu^*<\mu$
$$\theta\leq \sup\{ \sup \pcf_{\mu^*-\complete}(\frak{a}) \mid
\frak{a}\subseteq \Reg \cap (\mu^+,\chi) \text{ and } |\frak{a}|<\mu \}.$$
Just for every $\mu^*<\mu$ pick some $\chi_{\mu^*}$ such that
$$\theta\leq \sup\{ \sup \pcf_{\mu^*-\complete}(\frak{a}) \mid
\frak{a}\subseteq \Reg \cap (\mu^+,\chi_{\mu^*}) \text{ and } |\frak{a}|<\mu \},$$
and set $\chi=\bigcup_{\mu^*<\mu}\chi_{\mu^*}$.

So let us assume that $\cof(\theta)=\cof(\mu)$. Denote this common cofinality by $\kappa$.
By the assumption of the theorem $\kappa>\aleph_0$.
\\
Let $\l \mu_i \mid i<\kappa \r$ be an increasing continuous sequence with limit $\mu$ such that
each $\mu_i$ is a strong limit cardinal.
Let  $\theta>\mu$ be singular cardinal of cofinality $\kappa$.
Fix an increasing continuous sequence $\l \theta_i \mid i<\kappa \r $ with limit $\theta$ such that $\theta_0>\mu$.
\\Suppose that there are no $\sigma_*<\mu$  which satisfies the conclusion of the theorem.
In particular, for every $i<\kappa$, $ \mu_i$ cannot serve as $\sigma_*$.
Hence there is $\chi_i<\theta$ such that
$$\theta= \sup\{ \sup \pcf_{\mu_i-\complete}(\frak{a}) \mid
\frak{a}\subseteq \Reg \cap (\mu^+,\chi_i) \text{ and } |\frak{a}|<\mu \}.$$
So, for each $j<\kappa$, there is $\frak{a}_{i,j}\subseteq \Reg\cap (\mu^+,\chi_i)$
of cardinality less than $\mu$ such that  $\pcf_{\mu_i-\complete}(\frak{a}_{i,j})\not \subseteq \theta_j$.

Set $\theta_\kappa:=\theta$. For every $i\leq \kappa$, we apply Theorem 2.1 of \cite{Sh:460} to $\mu$ and $\theta_i$.
There is $\sigma_i^*<\mu$ such that
$$\text{ if } \frak{a} \subseteq \Reg \cap (\mu^+,\theta_i) \text{ and }  |\frak{a} |<\mu \text{ then }
\pcf_{\sigma_i^*-\complete}(\frak{a} )\subseteq \theta_i.$$

Define now by induction a sequence $\l i(n) \mid n<\omega\r$ such that
\begin{enumerate}
  \item
  $i(n)<i(n+1)<\kappa$,
  \item
  $\sigma^*_\kappa<\mu_{i(0)}$,
  \item $\sigma^*_{i(n)}<\mu_{i(n+1)}$,
  \item $\chi_{i(n)}<\theta_{i(n+1)}$.

\end{enumerate}

Let $i(\omega)=\bigcup_{n<\omega} i(n)$. Then $i(\omega)<\kappa$, since $\kappa$ is a regular above $\aleph_0$.
So $\theta_{i(\omega)}<\theta$.
Now, for every $j<\kappa$ and $n<\omega$ the following holds:
$$\frak{a}_{i(n),j}\subseteq \Reg\cap (\mu^+, \chi_{i(n)})\subseteq \Reg\cap (\mu^+, \theta_{i(n+1)})
\subseteq \Reg\cap (\mu^+, \theta_{i(\omega)}) \text{ and }$$ $$\pcf_{\sigma_{i(n+1)}^*-\complete}(\frak{a}_{i(n),j} )\subseteq \theta_{i(n+1)}<\theta_{i(\omega)}.$$

Let $n<\omega$ and $j \in (i(\omega),\kappa)$. Then by the choice of $\frak{a}_{i(n),j}$ the following holds:
$$\frak{a}_{i(n),j}\subseteq \Reg\cap (\mu^+, \chi_{i(n)})\subseteq \Reg\cap (\mu^+, \theta_{i(n+1)}) \text {and }
\pcf_{\mu_{i(n)}-\complete}(\frak{a}_{i(n),j} )\not \subseteq \theta_{j}.$$

By the choice of $\sigma_{i(n+1)}^*$, we have
$$\pcf_{\sigma_{i(n+1)}^*-\complete}(\frak{a}_{i(n),j} ) \subseteq \theta_{i(n+1)}.$$

By  \ref{lem-x} there is  $\frak{b}_{i(n),j} \subseteq \pcf_{\sigma_{i(n+1)}^*-\complete}(\frak{a}_{i(n),j} )$
such that\\
$|\frak{b}_{i(n),j} |< \sigma_{i(n+1)}^*<\mu_{i(n+2)}<\mu_{i(\omega)}$ and $\pcf_{\mu_{i(n)}-\complete}(\frak{b}_{i(n),j} )\not\subseteq \theta_j$.\\
Obviously, $\frak{b}_{i(n),j} \subseteq \Reg\cap(\mu^+, \theta_{i(n+1)})$, since \\
$\pcf_{\sigma_{i(n+1)}^*-\complete}(\frak{a}_{i(n),j} )\subseteq  \theta_{i(n+1)}$.

Apply Theorem 2.1 of \cite{Sh:460} to $\mu_{i(\omega)}$ (recall that it is a strong limit) and $\theta_{i(\omega)}$.
So, there is $\sigma_*<\mu_{i(\omega)}$ such that
$$\text{ if } \frak{b} \subseteq \Reg \cap (\mu^+_{i(\omega)},\theta_{i(\omega)}) \text{ and }  |\frak{b} |<\mu_{i(\omega)} \text{ then }
\pcf_{\sigma_*-\complete}(\frak{b} )\subseteq \theta_{i(\omega)}.$$
Now take $n_*<\omega$ with $\mu_{i(n_*)}>\sigma_*$.
Then $\frak{b}_{i(n_*),j}\subseteq \Reg \cap (\mu^+_{i(\omega)},\theta_{i(\omega)}) \text{ and }  |\frak{b}_{i(n_*),j} |<\mu_{i(\omega)}$, but
$\pcf_{\mu_{i(n_*)}-\complete}(\frak{b}_{i(n_*),j} )\not\subseteq \theta_j>\theta_{i(\omega)}$.
Which is impossible. Contradiction.
\\
$\square$

\section{Applications. }

Let $\kappa$ be a measurable cardinal, $U$ be a
$\kappa$--complete non-principle
ultrafilter over $\kappa$ and let
 $j_U:V \to M\simeq {}^\kappa V/U$ be the corresponding elementary
embedding.
Denote $j_U$ further simply by $j$.

\begin{lemma}\label{lem2.1}
Let $\mu>2^\kappa$ be a singular cardinal of cofinality $\kappa$.
Then $j(\mu)\geq \pp_{\Gamma(\kappa)}(\mu)$.

\end{lemma}
\pr Let $\lambda<\pp_{\Gamma(\kappa)}^+(\mu)$
be a regular cardinal.
Then, by Theorem \ref{pcf-main},
there is an increasing sequence of regular cardinals $\l \lambda_i \mid i<\kappa\r $ converging to $\mu$
such that $\lambda=\tcf(\prod_{i<\kappa}\lambda_i, <_{J_\kappa^{\bd}})$.
The ultrafilter $U$ clearly extends the dual to $J_\kappa^{\bd}$.
Hence $[\l \lambda_i \mid i<\kappa\r ]_U$ represents an ordinal below $j(\mu)$ of cofinality $\lambda$.
Hence $j(\mu)>\lambda$ and we are done.
\\
$\square$

Let us denote for a singular cardinal $\mu$ of cofinality
$\kappa$ by $\mu^*$ the least singular $\xi \leq \mu$ of
cofinality $\kappa$ above $2^\kappa$ such that $\pp_{\Gamma(\kappa)}(\xi)\geq \mu$.
\\Then, by \cite{Sh:g}(Sh 355, 2.3(3), p.57),
$\pp_{\Gamma(\kappa)}(\mu)\leq^+ \pp_{\Gamma(\kappa)}(\mu^*)$.

\begin{lemma}\label{lem2.2}
Let $\mu>2^\kappa$ be a singular cardinal of cofinality $\kappa$.
Then $j(\mu)\geq \pp_{\Gamma(\kappa)}(\mu^*)$.

\end{lemma}
\pr By \ref{lem2.1},
$j(\mu^*)\geq \pp_{\Gamma(\kappa)}(\mu^*)$.
But $\mu^*\leq \mu$, hence $j(\mu^*)\leq j(\mu)$.
\\
$\square$

\begin{lemma}\label{lem2.3}
Let $\mu>2^\kappa$ be a singular cardinal of cofinality $\kappa$.
Let $\eta, \mu<\eta<j(\mu)$ be a regular cardinal.
Then $\eta\leq \pp_{\Gamma(\kappa)}(\mu^*)$.

\end{lemma}
\pr
\\Let $\eta,\mu<\eta<j(\mu)$ be a regular cardinal.
Let $f_\eta:\kappa \to \mu$ be a  function
which represents $\eta $ in $M$, i.e. $[f_\eta]_U=\eta$.
We can assume that $\rng(f_\eta)\subseteq \Reg \cap ((2^\kappa)^+, \mu)$, since $|j(2^\kappa)|=2^\kappa$ and so
$j(2^\kappa)<\mu<\eta$.
Set $\tau:=U$--limit of  $\rng(f_\eta)$.\footnote{ It is possible to force a situation where
such $\tau<\mu$. Start with a $\eta^{++}$--strong $\tau, \kappa<\tau<\mu$.
Use the extender based Magidor to blow up the power of $\tau$ to $\eta^+$ simultaneously
changing the cofinality of $\tau$ to $\kappa$. The forcing satisfies $\kappa^{++}$-c.c., so it will not effect
$\pp$ structure of cardinals different from $\tau$.}
Then $\tau>2^\kappa$.
\\Note that $\cof(\tau)=\kappa$. Otherwise, $f_\eta$ is just a constant function mod $U$. Let $\delta$ be the constant value. Then $\delta<j(\delta)=\eta$. By elementarity $\delta$ must be a regular cardinal.
But then $j''\delta$ is unbounded in $\eta$, which means that $\eta $ is a singular cardinal. Contradiction.
\\
Denote $f(\alpha)$ by $\tau_\alpha$, for every $\alpha<\kappa$.
Then each $\tau_\alpha$ is a regular cardinal in the interval $((2^\kappa)^+, \tau)$
and $\tau=\lim_U\l \tau_\alpha\mid \alpha<\kappa\r$.
We have $\eta=\tcf(\prod_{\alpha<\kappa} \tau_\alpha, <_U).$\\
Note that once $U$ is not normal we cannot claim that the function $\alpha \mapsto \tau_\alpha$ is one to one.
So there is a slight tension between the true cofinalities of the sequence $\l \tau_\alpha \mid \alpha<\kappa \r$
and of the set $\{\tau_\alpha \mid \alpha<\kappa\}$.
\\We will show in the next lemma (\ref{lem-set-sequence}) that this does not effect $\pp_{\Gamma(\kappa)}(\tau)$.
\\Namely, $\eta=\tcf(\prod_{\alpha<\kappa} \tau_\alpha, <_U)$ implies
$\pp_{\Gamma(\kappa)}(\tau)\geq \eta>\mu$.\footnote{ Actually, the original definition of $\pp$ (\cite{Sh:g}II,Definition 1.1, p.41)
involves sequences rather than sets.}
\\
Then, by the choice of $\mu^*$, we have $\mu^*\leq \tau$
By \cite{Sh:g}(Sh 355, 2.3(3), p.57), \quad
$\pp_{\Gamma(\kappa)}(\mu^*)\geq\pp_{\Gamma(\kappa)}(\tau)$.
\\
$\square$

\begin{lemma}\label{lem-set-sequence}\footnote{ A version of this lemma was suggested by Menachem Magidor.}
Let $\kappa$ be a regular cardinal and $\tau$ be a singular cardinal of cofinality $\kappa$.
Then $$\pp_{\Gamma(\kappa)}(\tau)=\sup\{\tcf(\prod_{\alpha<\kappa}\tau_\alpha,<_I) \mid
\l \tau_\alpha \mid \alpha<\kappa\r \text{ is a sequence of regular cardinals with}$$$$ \lim_I\l \tau_\alpha \mid \alpha<\kappa\r=\tau, I \text{ is a } \kappa \text{ complete ideal over } \kappa \text{ which extends }
J_\kappa^{\bd}\}.$$

\end{lemma}
\pr
Clearly,
$$\pp_{\Gamma(\kappa)}(\tau)\leq\sup\{\tcf(\prod_{\alpha<\kappa}\tau_\alpha,<_I) \mid
\l \tau_\alpha \mid \alpha<\kappa\r \text{ is a sequence of regular cardinals with}$$$$ \lim_I\l \tau_\alpha \mid \alpha<\kappa\r=\tau, I \text{ is a } \kappa \text{ complete ideal over } \kappa \text{ which extends }
J_\kappa^{\bd}\}.$$
Just if $\eta=
tcf((\prod \frak{a}, <_J))$, where $ \frak{a} \text{ is a set of } \kappa \text{ regular cardinals unbounded in } \tau,$
$J \text{ is a } \kappa-\text{complete ideal on } \frak{a} \text{ which includes } J_{\frak{a}}^{\bd}$.
Then we can view $\frak{a}$ as a $\kappa$--sequence.
\\Let us deal with the opposite direction.
Suppose that
$\eta=\tcf(\prod_{\alpha<\kappa} \tau_\alpha, <_I)$, where \\
$\l \tau_\alpha \mid \alpha<\kappa\r \text{ is a sequence of regular cardinals with} \lim_I\l \tau_\alpha \mid \alpha<\kappa\r=\tau,\\ I \text{ is a } \kappa \text{ complete ideal over } \kappa \text{ which extends }
J_\kappa^{\bd}.$
Without loss of generality we can assume that $\kappa<\tau_\alpha<\tau$, for every $\alpha<\kappa$.
Set $\frak{a}=\{\tau_\alpha \mid \alpha<\kappa\}$.
Define a projection $\pi:\kappa \to \frak{a}$ by setting $\pi(\alpha)=\tau_\alpha$.
Let $$J:=\{X \subseteq \frak{a} \mid \pi^{-1}{}"X \in I \}.$$
Then $J$ will be a $\kappa$--complete ideal on $\frak{a}$ which extends $J_\frak{a}^{\bd}$.

Let us argue that $\eta=\tcf(\prod \frak{a}, <_J)$.
\\
Fix a scale $\l f_i \mid i<\eta\r$ which witnesses $\eta=\tcf(\prod_{\alpha<\kappa} \tau_\alpha, <_I)$.
Define for a function $f \in \prod_{\alpha<\kappa} \tau_\alpha$ a function $\bar{f}\in\prod_{\alpha<\kappa} \tau_\alpha$ as follows:
$$\bar{f}(\alpha)=\sup\{f(\beta) \mid \tau_\beta=\tau_\alpha \}.$$
Note that for every $\alpha<\kappa$, \quad $\bar{f}(\alpha)<\tau_\alpha$, since $\tau_\alpha$ is a regular cardinal above $\kappa$.
\\Consider the sequence $\l \bar{f}_i \mid i<\kappa \r$.
It need not be a scale, since the sequence need not be $I$--increasing.
But this is easy to fix. Just note that
 for every $i<\eta$ there will be $i',i\leq i'<\eta$, such that
$$f_i\leq \bar{f}_i\leq_{I}\bar{f}_{i'}.$$
Just given $i<\eta$, find some $i',i\leq i'<\eta$, such that
$\bar{f}_i\leq_{I}{f}_{i'}$. Then $\bar{f}_i\leq_{I}{f}_{i'}\leq \bar{f}_{i'}$.
Now by induction it is easy to shrink the sequence $\l \bar{f}_i \mid i<\kappa \r$
and to obtain an $I$--increasing subsequence  $\l g_\xi \mid \xi<\eta\r$ which is a scale
in  $(\prod_{\alpha<\kappa} \tau_\alpha,<_I)$.
\\
For every $\xi<\eta$ define $h_\xi \in\prod\frak{a}$ as follows:
$$h_\xi(\rho)=g_\xi(\alpha), \text{ if } \rho=\tau_\alpha, \text{ for some (every) } \alpha<\kappa.$$
It is well defined since $g_\xi(\alpha)=g_\xi(\beta)$ once $\tau_\alpha=\tau_\beta$.
\\Let us argue that $\l h_\xi \mid \xi<\eta\r$ is a scale in $(\prod\frak{a},<_J)$.
\\Clearly, $\xi<\xi'$ implies $h_\xi<_J h_{\xi'}$, since $g_\xi<_I g_{\xi'}$.
\\Let $h \in \prod\frak{a}$. Consider $g \in \prod_{\alpha<\kappa} \tau_\alpha$ defined by setting
$g(\alpha)=h(\tau_\alpha)$. There is $\xi<\eta$ such that $g<_Ig_\xi$.
Then $h<_J h_\xi$, since
$$\pi^{-1}"\{\rho \in \frak{a} \mid h(\rho)<h_\xi(\rho)\}\supseteq \{\alpha<\kappa \mid g(\alpha)<g_\xi(\alpha)\}.$$
\\
$\square$

\begin{theorem}\label{lem2.4}
 Let $\mu>2^\kappa$ be a singular cardinal of cofinality
$\kappa$.
\\Then $\pp_{\Gamma(\kappa)}(\mu^*)\leq j(\mu)<\pp_{\Gamma(\kappa)}(\mu^*)^+$.

\end{theorem}
\pr Note that $j(\mu)$ is always singular. Just $\mu$ is a singular
cardinal, hence $j(\mu)$ is a singular in $M$ and so in $V$. Now the
conclusion follows by \ref{lem2.2},\ref{lem2.3}.
\\
$\square$

We can deduce now an affirmative answer to a question of D. Fremlin for cardinals of cofinality $\kappa$:\footnote{ Readers interested only in a full answer to Fremlin's question can jump after the corollary directly to \ref{thm0-10}.
The non-strict inequality in its conclusion suffices.}

\begin{corollary}\label{cor1}
Let $W$ be a non-principal $\kappa$--complete ultrafilter on $\kappa$
and $j_{W}:V\to M_W$ the corresponding elementary embedding.
Then for every $\mu$ of cofinality $\kappa$, $|j(\mu)|=|j_{W}(\mu)|$.

\end{corollary}
\pr Let $\mu$ be a cardinal of cofinality $\kappa$.
If $\mu<2^\kappa$,   then
$2^\kappa<j_W(\mu)<j_W(2^\kappa)<(2^\kappa)^+$, for any
non-principal $\kappa$--complete ultrafilter $W$ on $\kappa$.
\\If $\mu>2^\kappa$, then, by \ref{lem2.4},
$\pp_{\Gamma(\kappa)}(\mu^*)\leq j(\mu)<\pp_{\Gamma(\kappa)}(\mu^*)^+$.
But recall that $j$ was the elementary embedding of an arbitrary
non-principal $\kappa$--complete ultrafilter $U$ on $\kappa$ and the bounds do not depend on it.
Hence if $W$ is an other  non-principal $\kappa$--complete ultrafilter on $\kappa$,
then
$\pp_{\Gamma(\kappa)}(\mu^*)\leq j_W(\mu)<\pp_{\Gamma(\kappa)}(\mu^*)^+$.
\\
$\square$

\begin{corollary}\label{cor2}
For every $\mu$ of cofinality $\kappa$, $|j(\mu)|=|j(j(\mu))|$.

\end{corollary}
\pr It follows from \ref{cor1}.
Just take $W=U^2$ and note that $j(j(\mu))=j_{U^2}(\mu)$.
\\
$\square$

Our next tusk will be to show that the fist inequality is really a strict inequality.

\begin{lemma}\label{lem2.7}
Let $\mu>2^\kappa$ be a singular cardinal of
cofinality $\kappa$.
  Then $pp_{\Gamma(\kappa)}(\mu)\leq
(pp_{\Gamma(\kappa)}(\mu))^M$.\footnote{ $(pp_{\Gamma(\kappa)}(\mu))^M$ stands for $pp_{\Gamma(\kappa)}(\mu)$ as computed in $M$. Note that it is possible to have $(pp_{\Gamma(\kappa)}(\mu))^M>pp_{\Gamma(\kappa)}(\mu)$, just as
$(2^\kappa)^M>2^\kappa$.}

\end{lemma}
\pr Let $\eta, \mu<\eta < pp_{\Gamma(\kappa)}^+(\mu)$ be a regular
cardinal.
\\
By Theorem \ref{pcf-main}, there is an increasing converging to
$\mu$ sequence $\l \eta_i \mid i<\kappa\r$ of regular
cardinals such that
$$\eta=\tcf(\prod_{i<\kappa}\eta_i,<_{J_{\kappa}^{\bd}}).$$
Note that both $\l \eta_i \mid i<\kappa\r$ and
$J_{\kappa}^{\bd}$ are in $M$. Also ${}^\kappa M \subseteq M$,
hence each function of the witnessing scale is in $M$, however the
scale itself may be not in $M$. Still we can work inside $M$ and
define a scale recursively using functions from the $V$-scale.
\\Thus let $\l f_\tau\mid \tau<\eta\r$ be a scale mod $J_{\kappa}^{\bd}$
which witnesses $\eta=\tcf(\prod_{i<\kappa}\eta_i,<_{J_{\kappa}^{\bd}}).$
Work in $M$ and define recursively an increasing
mod $J_{\kappa}^{\bd}$ sequence of functions
  $\l g_\xi \mid \xi<\eta'\r$
in $\prod_{i<\kappa}\eta_i$ as far as possible.
\\We claim first that $\cof(\eta')=\eta$, as computed in $V$.
Thus if $\eta<\cof(\eta')$, then there will be $\tau^*<\eta$ such that $f_{\tau^*}\geq_{J_{\kappa}^{\bd}}g_\xi$,
for every $\xi<\eta'$, since for every $\xi<\eta'$ there is $\tau<\eta$ such that $f_{\tau}\geq_{J_{\kappa}^{\bd}}g_\xi$.
But having $f_{\tau^*}\geq_{J_{\kappa}^{\bd}}g_\xi$, for all $\xi<\eta'$, we can continue and define $g_{\eta'}$
to be $f_{\tau^*}$.
\\If $\eta>\cof(\eta')$, then again there will be $\tau^*<\eta$ such that $f_{\tau^*}\geq_{J_{\kappa}^{\bd}}g_\xi$,
for every $\xi<\eta'$, and again  we can continue and define $g_{\eta'}$
to be $f_{\tau^*}$.
\\So  $\cof(\eta')=\eta$.
Let $\l \eta'_\tau \mid \tau<\eta\r$ be a cofinal in $\eta'$ sequence (in $V$).
Now, for every $\tau<\eta$ there is $\tau', \tau\leq \tau'<\eta$ such that
$f_{\tau}\not\geq_{J_{\kappa}^{\bd}}g_{\tau'}$, since the sequence $\l g_\xi \mid \xi<\eta'\r$
is maximal.
Hence there is $A_\tau\subseteq \kappa, |A_\tau|=\kappa$ such that
$f_{\tau}\upr A_\tau<_{J_{\kappa}^{\bd}}g_{\eta'_{\tau'}}\upr A_\tau$.
But $\eta>\mu>2^\kappa$, hence there is $A^*\subseteq \kappa$ such that for $\eta$ many $\tau$'s
we have $A^*=A_\tau$. Then for every $\tau<\eta$ there is $\tau'', \tau\leq \tau''<\eta$ such that
$f_{\tau}\upr A^*<_{J_{\kappa}^{\bd}}g_{\eta'_{\tau''}}\upr A^*$.
\\It follows that the sequence $\l g_\xi\upr A^* \mid \xi<\eta'\r$
is a scale in $\tcf(\prod_{i\in A^*}\eta_i,<_{J_{A^*}^{\bd}}).$
Hence, in $M$, $\eta'<pp_{\Gamma(\kappa)}^+(\mu)$. But $\cof(\eta')=\eta$, hence, in $M$, $\eta\leq \eta'<pp_{\Gamma(\kappa)}^+(\mu)$.
\\
$\square$

\begin{lemma}\label{lem0-11}
 Let $\mu>2^\kappa$ be a singular cardinal of cofinality
$\kappa$ such that $\mu^*=\mu$.\\Then $j(\xi)<\mu$ for every
$\xi<\mu$.
\end{lemma}
\pr Suppose otherwise. Then there is $\xi<\mu$ such that $j(\xi)\geq\mu$.
Necessarily $\xi>2^\kappa$.
Let
$\eta$ be a regular cardinal $\xi\leq\eta<\mu$. Pick a function
$f_\eta:\kappa \to \xi$ which represents $\eta $ in $M$.
Without loss of generality we can assume that $\min(\rng(f_\eta))>2^\kappa$.
Let
$\delta_\eta\leq \xi$ be the $U$--limit of $\rng(f_\eta)$. Then
$\cof(\delta_\eta)=\kappa$ and $j(\delta_\eta)>\eta$. Also $\eta
\leq \pp_{\Gamma(\kappa)}(\delta_\eta)$, by the definition of
$\pp_{\Gamma(\kappa)}(\delta_\eta)$. By Lemma \ref{lem2.2}, we have
$j(\delta_\eta)\geq pp_{\Gamma(\kappa)}((\delta_\eta)^*)$, and by
\cite{Sh:g}(Sh 355, 2.3(3), p.57), $\pp_{\Gamma(\kappa)}(\delta_\eta)\leq
\pp_{\Gamma(\kappa)}((\delta_\eta)^*)$. Set $$\delta:=\min\{\delta_\eta \mid
\xi\leq\eta<\mu \text{ and } \eta \text{ is a regular cardinal
}\}.$$ Then $\pp_{\Gamma(\kappa)}(\delta) \geq \pp_{\Gamma(\kappa)}(\delta_\eta)$, for
every regular $\eta,\xi\leq \eta <\mu$. But
$\pp_{\Gamma(\kappa)}(\delta_\eta)\geq \eta$. Hence $\pp_{\Gamma(\kappa)}(\delta) \geq
\mu$ which is impossible since $\mu^*=\mu$. Contradiction.
\\
$\square$

\begin{lemma}\label{lem2-100}
Let $\mu>2^\kappa$ be a singular cardinal of cofinality
$\kappa$.\\Then $\pp_{\Gamma(\kappa)}(\mu^*)<j(\mu)$.

\end{lemma}
\pr By \ref{lem2.2} we have $j(\mu)\geq \pp_{\Gamma(\kappa)}(\mu^*).$
\\Suppose that
 $j(\mu)=\pp_{\Gamma(\kappa)}(\mu^*)$.
 Then $\mu=\mu^*$, since by \ref{lem2.2} we have $j(\mu^*)\geq \pp_{\Gamma(\kappa)}(\mu^*).$
 By Theorem \ref{thm2}, there is $\sigma_*<\kappa$ such that
  $$\forall \chi<\mu (\mu>\sup\{\sup \pcf_{\sigma_*-\complete}(\frak{a})
 \mid \frak{a} \subseteq \Reg \cap (\kappa^+,\chi) \wedge |\frak{a}|<\kappa \}).$$
  Then, by elementarity,
  $$M\models
 \forall \chi<j(\mu) (j(\mu)>\sup\{\sup \pcf_{j(\sigma_*)-\complete}(\frak{a})
 \mid \frak{a} \subseteq \Reg \cap (j(\kappa^+),\chi) \wedge |\frak{a}|<j(\kappa) \}).$$
 Clearly, $j(\sigma_*)=\sigma_*$.
 Take $\chi=\mu$. Let $\eta$ be a regular cardinal (i.e. of $V$) such that
  $$(*)\quad  M \models j(\mu)>\eta >\sup\{\sup \pcf_{\sigma_*-\complete}(\frak{a})
 \mid \frak{a} \subseteq \Reg \cap (j(\kappa^+),\mu) \wedge |\frak{a}|<j(\kappa) \}.$$
 Note that there are such $\eta$'s since $j(\mu)$ is a singular cardinal of cofinality $\cof(j(\kappa))$.
 By Lemma \ref{lem2.3}, then $\eta \leq  \pp_{\Gamma(\kappa)}(\mu)$.
 Now, by Lemma \ref{lem2.7}, $\pp_{\Gamma(\kappa)}(\mu)\leq
(\pp_{\Gamma(\kappa)}(\mu))^M$. Hence
$M \models \eta\leq \pp_{\Gamma(\kappa)}(\mu)$.
But then there is $\frak{a} \in M$
 such that
$$M\models  \frak{a} \subseteq \Reg \cap (j(\kappa^+),\mu) \wedge |\frak{a}|=\kappa \wedge
\eta\leq \max \pcf_{\kappa-\complete}(\frak{a}).$$
 Which clearly contradicts $(*)$.
 \\
 $\square$

 So we proved the following:

\begin{theorem}\label{thm2.123}
 Let $\mu>2^\kappa$ be a singular cardinal of cofinality
$\kappa$.
\\Then $\pp_{\Gamma(\kappa)}(\mu^*)< j(\mu)<\pp_{\Gamma(\kappa)}(\mu^*)^+$.

\end{theorem}

Deal now with cardinals of arbitrary cofinality.

\begin{theorem}\label{thm0-10}
Let $\tau$ be a cardinal.
 Then either
\begin{enumerate}
\item $\tau<\kappa$ and then $j(\tau)=\tau$,\\ or
\item $\kappa\leq \tau \leq 2^\kappa$ and then $j(\tau)>\tau$,
 $2^\kappa<j(\tau)<(2^\kappa)^+$,\\ or
 \item $\tau\geq (2^\kappa)^+$ and then $j(\tau)>\tau$ iff
there is a singular cardinal ${\mu}\leq \tau$  of cofinality $\kappa$ above $2^\kappa$ such that
$\pp_{\Gamma(\kappa)}({\mu})\geq\tau$, and if $\tau^*$ denotes the least such $\mu$, then\\
$\tau\leq\pp_{\Gamma(\kappa)}(\tau^*)<j(\tau)<\pp_{\Gamma(\kappa)}(\tau^*)^+$.
\end{enumerate}

\end{theorem}
\pr
Suppose otherwise. Let $\tau$ be the least cardinal witnessing this.
Clearly then $\tau
>(2^\kappa)^+$. If $\cof(\tau)=\kappa$, then
we apply \ref{thm2.123} to derive the contradiction. Suppose that $\cof(\tau)\not =\kappa$.

\begin{claim}\label{cl0-10000}
There is a singular cardinal $\xi$ of cofinality $\kappa$
such that $j(\xi)>\tau$.

\end{claim}
\pr
Thus let $f_\tau:\kappa \to \tau$ be a
function which represents $\tau$ in $M$.
Without loss of generality we can assume that $$\nu\in \rng(f_\tau)\Rightarrow (\nu>2^\kappa \text{ and } \nu
\text{ is a cardinal }).$$
 Then
either $f_\tau$ is a constant function mod $U$ or
$\xi:=U$--limit
$\rng(f_\tau)$ has cofinality $\kappa$.
\\Suppose first that $f_\tau$ is a constant function mod $U$ with value $\xi$.
If $\xi=\tau$, then $j(\tau)=\tau$. Suppose that $\xi<\tau$. Then $j(\xi)=\tau>\xi$ and also $\xi$ is a cardinal above
$2^\kappa$.
By minimality of $\tau$ then $\xi^*$ exists and
$$\pp_{\Gamma(\kappa)}(\xi^*)<\tau=j(\xi)<\pp_{\Gamma(\kappa)}(\xi^*)^+.$$
But this is impossible since $\tau$ is a cardinal.
Contradiction.
So  $\cof(\xi)=\kappa$ and $j(\xi)>\tau$.
\\
$\square$ of the claim.

Let  $\mu\leq \tau$ be
the least singular cardinal above $2^\kappa$ of cofinality $\kappa$ such that
$j(\mu)>\tau$. We claim that $\mu=\mu^*$. Note that by \ref{thm2.123},
we have $\pp_{\Gamma(\kappa)}(\mu^*)<j(\mu^*)\leq j(\mu)<\pp_{\Gamma(\kappa)}(\mu^*)^+$.
 \\$\tau$ is a cardinal below $j(\mu)$,  hence $\tau\leq
\pp_{\Gamma(\kappa)}(\mu^*)<j(\mu^*)$. The minimality of $\mu$ implies then that $\mu=\mu^*$.
Note that also $\tau^*=\mu$. Thus $\pp_{\Gamma(\kappa)}(\tau^*)\geq \tau\geq \mu=\mu^*$, and so
$\tau^*\geq\mu$. Also $\tau\leq
\pp_{\Gamma(\kappa)}(\mu)$ implies $\tau^*\leq\mu$.
\\
 Apply finally
\ref{cor2}. It follows that $|j(j(\mu))|=|j(\mu)|$, but
$j(\mu)>\tau$, hence $j(j(\mu))>j(\tau)>j(\mu)$.
So
$$\pp_{\Gamma(\kappa)}({\mu})<j(\mu)<j(\tau)<\pp_{\Gamma(\kappa)}({\mu})^+,$$
 and we are done.
\\
$\square$

Now affirmative answers to a question of D. Fremlin and to questions 4,5 of \cite{G-Sh}
follow easily.\footnote{ Non strict inequality $\pp_{\Gamma(\kappa)}(\tau^*)\leq j(\tau)<\pp_{\Gamma(\kappa)}(\tau^*)^+$
suffices for a question of D. Fremlin and 4 of \cite{G-Sh}. }
\begin{corollary}\label{cor3}
Let $W$ be a non-principal $\kappa$--complete ultrafilter on $\kappa$
and $j_{W}:V\to M_W$ the corresponding elementary embedding.
Then for every $\tau$, $|j(\tau)|=|j_{W}(\tau)|$.

\end{corollary}
\pr
Let $W$ be a non-principal $\kappa$--complete ultrafilter on $\kappa$
and $j_{W}:V\to M_W$ the corresponding elementary embedding.
Let $\tau$ be an ordinal. Without loss of generality we can assume that $\tau$ is a cardinal,
otherwise just replace it by $|\tau|$.
Now by \ref{thm0-10},
 $j(\tau)>\tau$ iff $j_W(\tau)>\tau$
and if  $j(\tau)>\tau$ then either
 $j(\tau), j_W(\tau) \in (2^\kappa,(2^\kappa)^+)$,\\ or
 $j(\tau), j_W(\tau) \in (\pp_{\Gamma(\kappa)}(\tau^*),\pp_{\Gamma(\kappa)}(\tau^*)^+)$.
 \\
$\square$

\begin{corollary}\label{cor4}
For every $\tau$, $|j(\tau)|=|j(j(\tau))|$.

\end{corollary}
\pr
Apply \ref{cor3} with $W=U^2$.
\\
$\square$

It is straightforward to extend this to arbitrary iterated ultrapowers of $U$:

\begin{corollary}\label{cor4}
Let $\tau$ be a cardinal with $j(\tau)>\tau$.
Let $\alpha\leq 2^\kappa$, if $\tau\leq 2^\kappa$, and $\alpha\leq \pp_{\Gamma(\kappa)}(\tau^*)$, if
$\tau> 2^\kappa$. Then
$|j(\tau)|=|j_\alpha(\tau))|$,
where $j_\alpha:V\to M_\alpha$ denotes the $\alpha$-th iterated ultrapower of $U$.

\end{corollary}

\begin{corollary}\label{cor5}
For every $\tau$, if $j(\tau)\not = \tau$, then $j(\tau)$ is not a cardinal.

\end{corollary}
\pr Follows immediately from \ref{thm0-10}.
\\
$\square$

The following question looks natural:

\emph{Let $\alpha$ be any ordinal. Suppose $j(\alpha)>\alpha$.
Let $W$ be a non-principal $\kappa$--complete ultrafilter on $\kappa$
and $j_{W}:V\to M_W$ the corresponding elementary embedding.
Does then  $j_W(\alpha)>\alpha$?}

Next statement answers it negatively assuming that $o(\kappa)$-- the Mitchell order  of $\kappa$ is at least 2.

\begin{proposition}\label{prop200}
Let $W$ be a non-principal $\kappa$--complete ultrafilter on $\kappa$
and $j_{W}:V\to M_W$ the corresponding elementary embedding.
Suppose that  $U \triangleleft W $, i.e.  $U \in M_W$.
Then  $j_W(\alpha)>\alpha=j(\alpha)$,  for some $\alpha<(2^\kappa)^+$.

\end{proposition}
\pr
Let $\alpha=j_\omega(\kappa)$, i.e. the $\omega$-th iterate of $\kappa$ by $U$.
Then $j(\alpha)=\alpha$, since $j_\omega(\kappa)=\cup_{n<\omega}j_n(\kappa)$.
Let us argue that $j_W(\alpha)>\alpha$. Thus we have $U$ in $M_W$.
So $j_\omega(\kappa)$ as computed in $M_W$ is the real $j_\omega(\kappa)$.
In addition
$$M_W \models |j_\omega(\kappa)|=2^\kappa<(2^\kappa)^+<j_W(\kappa),$$
and so $\kappa<\alpha=j_\omega(\kappa)<j_W(\kappa)$.
Hence
$$j_W(\alpha)=j_W(j_\omega(\kappa))>j_W(\kappa)>\alpha.$$
$\square$

Let us note that the previous proposition is sharp.

\begin{proposition}\label{prop201}
Suppose that there is no inner model with a measurable of the Mitchell order $\geq 2$.
Let $W$ be a non-principal $\kappa$--complete ultrafilter on $\kappa$
and $j_{W}:V\to M_W$ the corresponding elementary embedding.
Then $j(\alpha)>\alpha$ iff $j_W(\alpha)>\alpha$,  for every ordinal $\alpha$.

\end{proposition}
\pr
Assume that $U$ is normal or just replace it by such.
Let $W$ be a non-principal $\kappa$--complete ultrafilter on $\kappa$
and $j_{W}:V\to M_W$ the corresponding elementary embedding.
\\The assumption that there is no inner model with a measurable of the Mitchell order $\geq 2$
guarantees that there exists the core model. Denote  denote it by $\calK$. Let $U^*= U \cap \calK$.
Then it is a normal ultrafilter over $\kappa$ in $\calK$.
Denote by $j^*$ its elementary embedding.
 Then $j_W \upr \calK=j_n^*$, for some $n<\omega$, since ${}^\omega M_W \subset M_W$ there are no measurable
 cardinals in $\calK$ of the Mitchell order 2.
\\Hence we need to argue that
$$j^*(\alpha)>\alpha\Leftrightarrow j_n^*(\alpha)>\alpha,$$  for every ordinal $\alpha$ and every $n<\omega$.
But this is trivial, since $j^*(\alpha)>\alpha$ implies
$j_2^*(\alpha)=j^*(j^*(\alpha))>j^*(\alpha)>\alpha$ and in general
$j_{k+1}^*(\alpha)=j^*(j_k^*(\alpha))>j_k^*(\alpha)>\alpha$, for every $k,0<k<\omega$.
On the other hand, if $j^*(\alpha)=\alpha$, then $j_\xi^*(\alpha)=\alpha$, for every $\xi$.
\\
$\square$

\section{ Concluding remarks and open problems.}

\textbf{Question 1.} \emph{Is weak compactness really needed for Theorem
\ref{pcf-main}? Or explicitly:
\\Let $\kappa$ a regular cardinal.
Let $\mu>2^\kappa$ be a singular cardinal of cofinality $\kappa$.
Suppose that $\lambda<\pp^+_{\Gamma(\kappa)}(\mu)$. Is there an
increasing sequence $\l \lambda_i \mid i<\kappa \r$ of regular
cardinals converging to $\mu$ such that $\lambda=
\tcf(\prod_{i<\kappa} \lambda_i, <_{J^{bd}_\kappa})$?}

See \cite{Sh:g} pp.443-444, 5.7 about the related results.
\\
\textbf{Question 2.} \emph{ Does Theorem \ref{thm2} remain true assuming $\cof(\mu)=\cof(\theta)=\omega$?}

Suppose now that we have an $\omega_1$-saturated $\kappa$--complete ideal on $\kappa$ instead of a
$\kappa$--complete ultrafilter. The following generic analogs of questions 4,5 of \cite{G-Sh} and of a question
of Fremlin are natural:
\\
\textbf{Question 3.}\emph{ Let $W$ be an $\omega_1$--saturated filter on
$\kappa$. Does each the following hold:}

\begin{enumerate}
  \item $\Vdash_{W^+} \forall \tau (\lusim{j}_W(\tau)>\tau\longrightarrow \tau \text{ \emph{is not a cardinal}). }$
  \item $\Vdash_{W^+} \forall \tau (|\lusim{j}_W(\tau)|=|\lusim{j}_W(\lusim{j}_W(\tau))|). $
  \item\emph{ Let $W_1$ be an other $\omega_1$--saturated filter on
$\kappa$.
Suppose that for some $\tau$ we have $\delta,\delta_1$ such that}
\begin{itemize}
  \item $\Vdash_{W^+} \lusim{j}_W(\tau)=\check{\delta}$,
  \item $\Vdash_{W^+_1} \lusim{j}_{W_1}(\tau)=\check{\delta}_1$.
\end{itemize}
\emph{Then $|\delta|=|\delta_1|$.}
\end{enumerate}

Note that in such situation $2^{\aleph_0}\geq \kappa$ and so \ref{pcf-main} does not apply.
Assuming   variations of SWH and basing on \cite{Sh:g}, Sh371,
it is possible to answer positively this questions for $\tau>2^\kappa$.

Recall a question of similar flavor from \cite{G-Sh} (Problem 6):
\\
\textbf{Question 4.} \emph{Let $W$ be an $\omega_1$--saturated filter on
$\kappa$. Can the following happen:}

 $\Vdash_{W^+} \lusim{j}_W(\kappa)$ \emph{is a cardinal?
Or even }
$\Vdash_{W^+} \lusim{j}_W(\kappa)=\kappa^{++}$?

\end{document}